\documentclass{article}

\usepackage{arxiv}

\usepackage{dsfont}
\usepackage{amssymb}
\usepackage{amsmath}
\usepackage{amsfonts}
\usepackage{mathtools}
\usepackage{empheq}
\usepackage[utf8]{inputenc} % allow utf-8 input
\usepackage[T1]{fontenc}    % use 8-bit T1 fonts
\usepackage{hyperref}       % hyperlinks
\usepackage{url}            % simple URL typesetting
\usepackage{booktabs}       % professional-quality tables
\usepackage{amsfonts}       % blackboard math symbols
\usepackage{nicefrac}       % compact symbols for 1/2, etc.
\usepackage{microtype}      % microtypography
\usepackage{cleveref}       % smart cross-referencing
\usepackage{lipsum}         % Can be removed after putting your text content
\usepackage{graphicx}
\usepackage{doi}
\usepackage{yhmath}

\newtheorem{lemma}{Lemma}

\newtheorem{theorem}{Theorem}
\newtheorem{definition}{Definition}
\newtheorem{remark}{Remark}

\title{Invariant cones for strange attractors of Lozi, H\'{e}non and Belykh type maps}

\author{D.~A.~Grechko \\
	Volga State University of Water Transport\\
	Nizhny Novgorod, 603950 Russia \\
        Lobachevsky State University of Nizhny Novgorod \\
        Nizhny Novgorod, 603022 Russia \\
	\texttt{d.grechko.18@gmail.com} \\
	\And
V.~N.~Belykh \\
	Volga State University of Water Transport\\
	Nizhny Novgorod, 603950 Russia \\
        Lobachevsky State University of Nizhny Novgorod \\
        Nizhny Novgorod, 603022 Russia \\
	\texttt{Belykh.vn@vsuwt.ru} \\
	\AND
	N.~V.~Barabash \\
	Volga State University of Water Transport\\
	Nizhny Novgorod, 603950 Russia \\
        Lobachevsky State University of Nizhny Novgorod \\
        Nizhny Novgorod, 603022 Russia \\
	\texttt{barabash@itmm.unn.ru} \\
}

\renewcommand{\shorttitle}{Invariant cones for strange attractors of Lozi, H\'{e}non and Belykh type maps}

\hypersetup{
pdftitle={Invariant cones for strange attractors of Lozi, H\'{e}non and Belykh type maps},
pdfsubject={math.DS},
pdfauthor={D.~A.~Grechko, V.~N.~Belykh,N.~V.~Barabash},
pdfkeywords={dynamical systems, nonlinear maps, attractors, hyperbolicity, bifurcations},
}
\begin{document}
\maketitle

\begin{abstract}
We consider strange attractors of two dimensional generalized map with one nonlinearity \cite{BelykhMatSbornik1995,BelykhDinamicaSystem1976} such that Lozi, H\'{e}non and Belykh maps are particular cases of it. We describe technique of invariant expanding and contracting cones creation for study of hyperbolic properties. Theorems of singular hyperbolic attractors for new modifications of Lozi, H\'{e}non and Belykh-type maps are presented.
\end{abstract}

% keywords can be removed
\keywords{dynamical systems \and nonlinear maps \and attractors \and hyperbolicity \and bifurcations}

\section{Introduction}

In this paper we consider map
$f:\;\mathbb{R}^2\rightarrow\mathbb{R}^2$ of the form \cite{BelykhMatSbornik1995, BelykhDinamicaSystem1976}
\begin{equation}
f:\quad
    \begin{array}{l}
    \bar{x}=x+y+g(x),\\
    \bar{y}=\lambda\big(y+g(x)\big),
    \end{array}
    \label{sys:map_f}
\end{equation}
where $0<\lambda<1$ and $g(x)$ is a piecewise smooth function.
This map arised as a model of phase-locked loop (PLL) \cite{BelykhDinamicaSystem1976,PLLoverview2007}.
Namely, the simplest version of PLL models has the form of ODE system
\begin{equation}
    \begin{array}{l}
         \dot{x}=u+\omega_1g(x),\\
         \dot{u}=-\lambda_0u+\omega_2g(x), 
    \end{array}
    \label{sys:PLL}
\end{equation}
where $\omega_{1,2}$, $\lambda_0$ are positive paramters, and $g(x)$ is sine-like periodic function serving the output of PLL phase detector.

The discrete time model of PLL is similar to result of applying Euler method with discretely defined time derivatives in \eqref{sys:PLL}
\begin{equation}
    \begin{array}{l}
         x_{i+1}=x_i+hu_i+h\omega_1g(x_i),\\
         u_{i+1}=u_i-h\lambda_0u_i+h\omega_2g(x_i).
    \end{array}
    \label{sys:PLL_map}
\end{equation}
The variable $u$ and parameters $\lambda_0$, $\omega_1$, $\omega_2$ change in \eqref{sys:PLL_map}
\begin{equation}
    hu=y,\quad \lambda=1-\lambda_0h,\quad
    h\omega_1=a,\quad \omega_2h^2=\lambda a
    \label{eq:change_for_Belykn_map}
\end{equation}
leads to map~\eqref{sys:map_f}.

Now we demonstrate that this map has a universal character. 
\subsection{Sine-Gordon equation}
For the steady-state equation of Sine-Gordon type
\begin{equation*}
    (1-\lambda)\dfrac{\partial u}{\partial x}
    +
    \dfrac{\partial^2 u}{\partial x^2}=\sin u
\end{equation*}
we consider the following discrete version
\begin{equation*}
    (1-\lambda)(u_j-u_{j-1})+u_{j+1}-2u_{j}+u_{j-1}=\sin u_j
\end{equation*}
or equation
\begin{equation}
    u_{j+1}-u_j-\lambda(u_j-u_{j-1})=\sin u_j.
\label{sys:map_sine_gordon}
\end{equation}
If we set $u_j=x_j$, $\lambda(u_j-u_{j-1})=y_j$
then equation~\eqref{sys:map_sine_gordon} implies that
\begin{equation*}
\begin{array}{l}
     x_{j+1}=x_j+y_j+\sin x_j,\\
     y_{j+1}=\lambda(y_j+\sin x_j), 
\end{array}
\end{equation*}
and we obtain the dynamical system generated by map~\eqref{sys:map_f} with $g(x)=\sin x$.

\subsection{Standard map}
For $\lambda=1$, $g(x)=k\sin x$ map~\eqref{sys:map_f} becomes area preserving Chirikov standard map \cite{Chirikov1979} playing essential role in theory of the conservative chaos.

\subsection{Zaslavsky map}
Zaslavsky map \cite{Zaslavsky1978,SagdeevAfraim1990} has the form 
\begin{equation*}
    \begin{array}{l}
         \bar{z}=\eta(z+a\sin\theta),\\
         \bar{\theta}=\theta+\omega+z+a\sin\theta,
    \end{array}
\end{equation*}
where $a$, $\eta$ and $\omega$ are parameters.
If we set $\theta=x$, $z=y-\omega\eta$, $g(x)=\omega(1-\eta)+a\sin x$, $\eta=\lambda$
then Zaslavsky map takes the form~\eqref{sys:map_f}.

\subsection{Lozi-H\'{e}non map}
Consider famous 2-D map 
\begin{equation}
    \begin{array}{l}
         \bar{x}=y+1-aU(x),\\
         \bar{y}=bx, 
    \end{array}
    \label{sys:Henon-Lozi}
\end{equation}
which is H\'{e}non map \cite{Henon1976,CarlesonBenedicks1991} for $U(x)=x^2$ and the Lozi map \cite{Lozi1978} for $U(x)=|x|$.
Variable and parameter change in \eqref{sys:Henon-Lozi}
\begin{equation*}
    (x,y)\rightarrow(x, y-bx),\quad \lambda=-b,
\end{equation*}
leads the Lozi-H\'{e}non map to the form \eqref{sys:map_f} with the following function
\begin{equation*}
    g(x)=1-(1+\lambda)x-aU(x).
\end{equation*}
In this paper which is based on the results from \cite{BelykhMatSbornik1995} we present sufficient conditions of hyperbolicity for map~\eqref{sys:map_f} (section~\ref{sec:hyper}).
As an example of applying these conditions we consider the Lozi map written in the form of map~\eqref{sys:map_f} (section~\ref{sec:Lozi}).
In section~\ref{sec:Lozi_Hybrid} we present analysis of nonlinear piecewise smooth hybrid Lozi-H\'{e}non map.
Finally in section~\ref{sec:Belykh} we consider the original form of Belykh map and present new features of Belykh attractor in the case of periodic function $g(x)$.

\section{\label{sec:hyper} Conditions of hyperbolicity}

Invariant expanding and contracting cones play an essential role in the theory of hyperbilocity \cite{AnosovGeodes1967,AnosovGrinesDinamSys1991,KatokHasselblatt1999,AfraimovichChernovSataev1995}, being the main tool of proven hyperbolic properties of certain maps \cite{Glindinning2021,Misiurewicz1980,BelykhMatSbornik1995,BelykhKomrakovUkrainsky2001,BelykhHybridSys2010,BelykhGrechkoDinSys2018,BarabashEuropJ2020}. 
In this section we present detailed technique of cones construction briefly presented in \cite{BelykhMatSbornik1995}.

\begin{definition}
Assume that a smooth map $F:\;G\rightarrow G$, where $G$ is an open subset of $\mathbb{R}^2$, has the same in each point of $G$ invariant cones $K^u$ and $K^s$ \cite{BelykhMatSbornik1995} such that $DFK^u\subset K^u$ and $DF^{-1}K^s\subset K^s$, where linearization $DF$ is expanding in $K^u$ and contracting in $K^s$ for any point of $G$.
Then map $F$ is called hyperbolic.
\end{definition}

Note that despite of usual notation of hyperbolicity \cite{KatokHasselblatt1999,AfraimovichLecture2002} this definition has a strong restriction that cones $K^u$ and $K^s$ are one and the same in each point of region $G$.

Consider map~\eqref{sys:map_f} in
region $G=G^-\cup G^+$ where $G^{-(+)}$ corresponds to pieces of function $g(x)$ monotonicity with negative (positive, respectively) derivative $g'(x)$.
In other words, we consider piecewise smooth map~\eqref{sys:map_f} with bounded away from zero derivative $g'(x)$
\begin{equation*}
g'(x)\triangleq d(x), \quad |d(x)|\geq\varepsilon>0,
\end{equation*}
for each point of $G$. 

\begin{theorem}
    Let the conditions
    \begin{equation}
    \begin{array}{ll}
         d(x)\geq \varepsilon>0 &\textrm{for } x\in G^+,\\
         \\
         d(x)<-2(1+\lambda) &\textrm{for } x\in G^-, 
    \end{array}
    \label{eq:theorem_cond_d}
    \end{equation}
    hold. 
    Then map~\eqref{sys:map_f} is hyperbolic in region $G$ covered by the next invariant unstable and stable cones
    \begin{equation}
    K^u=\left\lbrace 0< \frac{u_2}{u_1}<\lambda \right\rbrace,
    \quad
        K^s=\left\lbrace -\frac{1}{1-\lambda}< \frac{u_1}{u_2}<\frac{1}{1+\lambda}\right\rbrace,
        \label{eq:cones}
\end{equation}
where $u_{1,2}$ are local coordinates in each point of $G$.
    \label{th:Hyper}
\end{theorem}

\textit{Proof.}
To construct cones $K^{u,s}$ we consider variational equations as a linearized map
\begin{equation}
    L:\quad\bar{\textbf{u}}=A(x)\textbf{u},
    \label{sys:map_f_L}
\end{equation}
where
\begin{equation}
\textbf{u}=
\begin{pmatrix}
u_1\\
u_2
\end{pmatrix},\quad
\bar{\textbf{u}}=\begin{pmatrix}
\bar{u}_1\\
\bar{u}_2
\end{pmatrix},
\quad
A(x)=\begin{pmatrix}
1+d(x)& 1 \\
\lambda d(x) &\lambda
\end{pmatrix}.
\label{eq:local_u_12}
\end{equation}

%Note, that map~\eqref{sys:map_f} is dissipative due to $\det A(x)=\lambda<1$. 
Eigenvalues $\mu_{1,2}(x)$ and eigenvectors $\textbf{V}_{1,2}(x)=\big(1, \alpha_{1,2}(x)\big)$ of matrix $A(x)$ take the form
\begin{equation}
\begin{array}{l}
 \mu_{1,2}(x)=\frac{1}{2}\left(1+d(x)+\lambda\pm\sqrt{(1+d(x)+\lambda)^2-4\lambda}\right),\\
 \\
\alpha_{1,2}(x)=\frac{1}{2}\left( -1-d(x)+\lambda\pm\sqrt{(1+d(x)+\lambda)^2-4\lambda}\right). 
\end{array}
\label{eq:multypliers_vectors}
\end{equation}

Under conditions~\eqref{eq:theorem_cond_d} of Theorem~\ref{th:Hyper} eigenvalues \eqref{eq:multypliers_vectors} satisfy
the following conditions
\begin{equation}
    \begin{array}{ccl}
          \mu_1(x)>1,& 0<\mu_2(x)<1,& \textrm{for}\;d(x)>0,\\
          \\
    -1<\mu_1(x)<0, &\mu_2(x)<-1,& \textrm{for}\;d(x)<-2(1+\lambda),
    \end{array}
    \label{eq:conditions_d>0}
\end{equation}
which imply that each point of $G$ is of saddle type.

From formulas \eqref{eq:multypliers_vectors} satisfying conditions~\eqref{eq:conditions_d>0} it follows that the eigenvectors define the next directions
\begin{equation}
    \begin{array}{ll}
    \textbf{V}_1\big\vert_{G^+}\;\textrm{and}\;\textbf{V}_2\big\vert_{G^-}
    \;\textrm{--  unstable directions}, \\
    \\
        \textbf{V}_2\big\vert_{G^+}\;\textrm{and}\;\textbf{V}_1\big\vert_{G^-}
    \;\textrm{--  stable directions}.
    \end{array}
    \label{eq:V_properties}
\end{equation}

The angle boundaries of eigenvectors are defined by edge values of coordinates $\alpha_{1,2}(x)$.
Directions of the angle turn are defined by derivatives of $\alpha_{1,2}(x)$ with respect to $d(x)$
\begin{equation*}
   \alpha_{1,2}'(x)=
   \frac{1}{2}\left( -1\pm\frac{1+d(x)+\lambda}{\sqrt{(1+d(x)+\lambda)^2-4\lambda}} \right)
\end{equation*}
which gives inequalities
\begin{equation*}
    \begin{array}{rrl}
    \alpha_1'(x)>0,& \alpha_2'(x)<0&\textrm{for}\;d\geq\varepsilon>0,\\
    \\
     \alpha_1'(x)<0,& \alpha_2'(x)>0&\textrm{for}\;d<-2(1+\lambda).
    \end{array}
\end{equation*}

Due to \eqref{eq:multypliers_vectors} we obtain the next intervals of monotonic change of $\alpha_{1,2}$ as functions of $d(x)$
\begin{equation}
\begin{array}{lll}
\begin{array}{l}
\alpha_1(x)\in(0,\lambda),\\
 \alpha_2(x)\in(-1+\lambda,-\infty), 
\end{array}
& \textrm{for}\;d(x)\in(0,\infty),\\
\\
\begin{array}{l}
 \alpha_1(x)\in(\infty,1+\lambda),\\
 \alpha_2(x)\in(\lambda,2\lambda),
\end{array}
&\textrm{for}\;d(x)\in(-\infty,-2(1+\lambda))
\end{array}
\label{eq:alfa_bound}
\end{equation}

\begin{figure}
    \centering
    \includegraphics[width=0.75\textwidth]{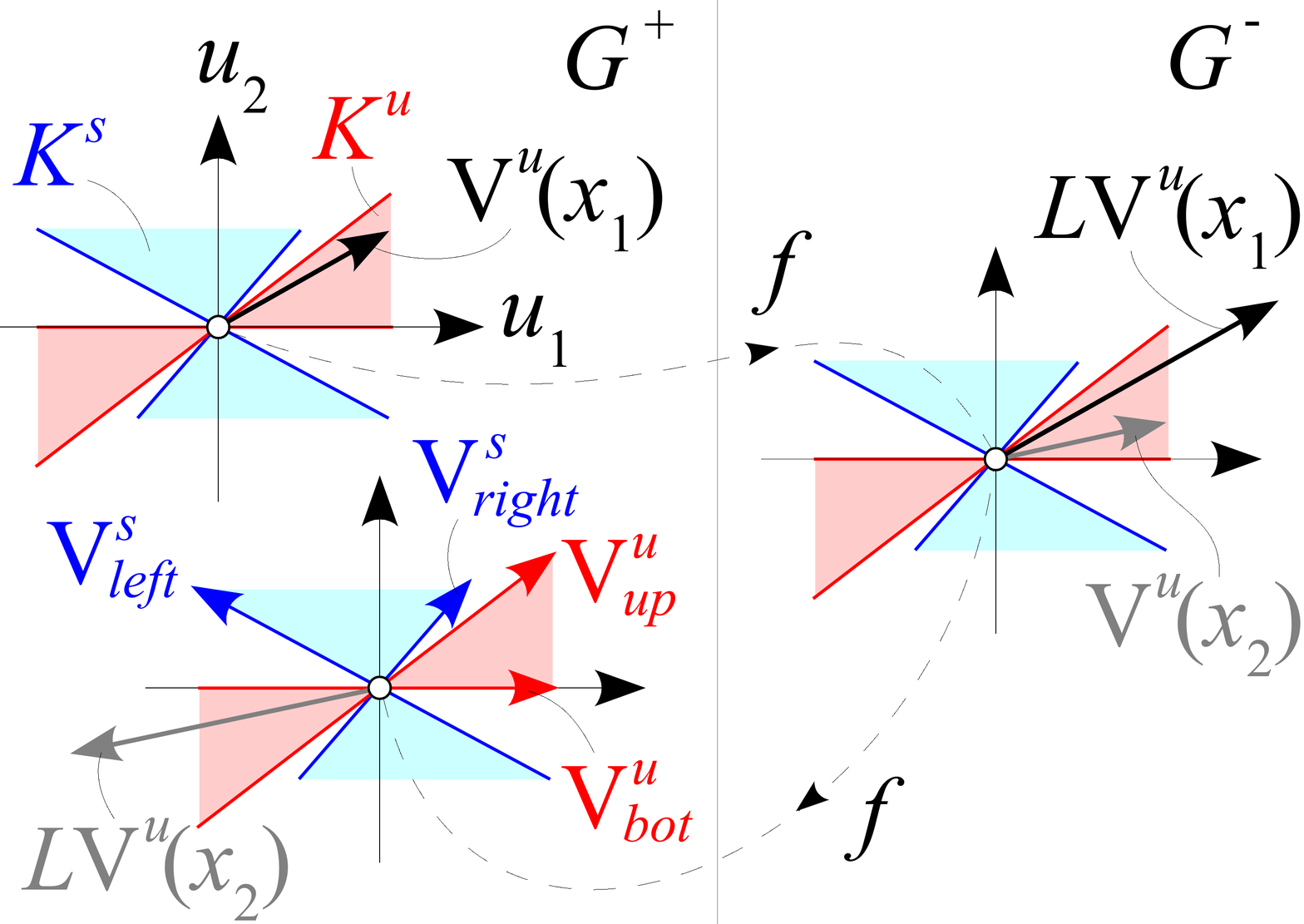}
    \caption{Unstable $K^u$ and stable $K^s$ invariant cones (red and blue triangles, respectively).
    Edge vectors of cones are $\textbf{V}^u_{up}$, $\textbf{V}^u_{bot}$ (red arrows) for  $K^u$, and $\textbf{V}^s_{left}$, $\textbf{V}^s_{right}$ (blue arrows) for  $K^s$.
    In depicted case, while point $(x_1,y_1)$ is mapped by $f$ (dashed black curves) from $G^+$ to $G^-$, its eigenvector $\textbf{V}^u(x_1)\in K^u$ (black bold arrow) is expanded under linear map $L$ (see Eq.~\eqref{sys:map_f_L}).
    In the next step, point $(x_2,y_2)$ returns to $G^+$ under the action of $f$, and its unstable eigenvector $\textbf{V}^u(x_2)\in K^u$ (gray bold arrow) expands under $L$ rotated by $\pi$ (due to negative multipliers, see Eq.~\eqref{eq:multypliers_vectors}). 
    In stable cone $K^s$ all eigenvectors are mapped with contraction (not shown).}
    \label{fig:acting_cones}
\end{figure}

Thus, the family of all possible eigenvectors $\textbf{V}^u$ and $\textbf{V}^s$ in region $G$ for $d(x)>0$ and $d(x)<-2(1+\lambda)$
are bounded by the next edge vectors (see Fig.~\ref{fig:acting_cones})
\begin{equation}
\begin{array}{ll}
\textbf{V}_{bot}^u=(1, 0),&
\textbf{V}_{up}^u=(1, 2\lambda),\\
\\
\textbf{V}_{left}^s=(-1, 1-\lambda),&
\textbf{V}_{right}^s=(1, 1+\lambda).
\end{array}
\end{equation}

As a result using these edge vectors we obtain the form of unstable and stable invariant cones \eqref{eq:cones}. 
$\Box$

Figure~\ref{fig:acting_cones} illustrates iterates of linearized map $L$ along a trajectory of $f$ (eq.~\eqref{sys:map_f}) in $G$.
Due to Theorem~\ref{th:Hyper} any eigenvector $\textbf{V}^u(x)$ \big($\textbf{V}^s(x)$\big) belongs to unstable (stable) invariant cone $K^u$ ($K^s$).
Therefore, any trajectory in region $G$ is of saddle type. 

Using Theorem~\ref{th:Hyper} we obtain 
the next statement.
\begin{theorem}
    Let $G$ be a trapping domain, $fG\subset G$, of map $f$ (eq.~\eqref{sys:map_f}), satisfying conditions of Theorem~\ref{th:Hyper}.
    Then map $f$ has singular hyperbolic attractor.
    \label{th:attractor}
\end{theorem}

\textit{Proof.}
As $G$ is trapping domain, so map~\eqref{sys:map_f} has an attractor $\mathcal{A}=\bigcap_{k=1}^\infty f^kG$.
Then from Theorem~\ref{th:Hyper} it follows that attractor $\mathcal{A}$ is hyperbolic.
The singularity of this attractor appears when points of attractor $\mathcal{A}$ lie at discontinuity lines $\left\lbrace x=\textrm{const}, y\in\mathbb{R}^1\right\rbrace$ at which derivative $d(x)$ changes the sign.
$\Box$

\section{\label{sec:Lozi}Continuous piecewise linear function and Lozi attractor}

We start to demonstrate application of Theorem~\ref{th:Hyper} for the Lozi map written in the form~\eqref{sys:map_f} with piecewise linear function
\begin{equation}
    g(x)=1-a|x|-(1+\lambda)x,\quad\lambda=-b,
     \label{eq:g(x)_Lozi}
\end{equation}
where $\lambda>0$ corresponds to orientation preserving case. 

\begin{figure}
    \centering
    \includegraphics[width=0.75\textwidth]{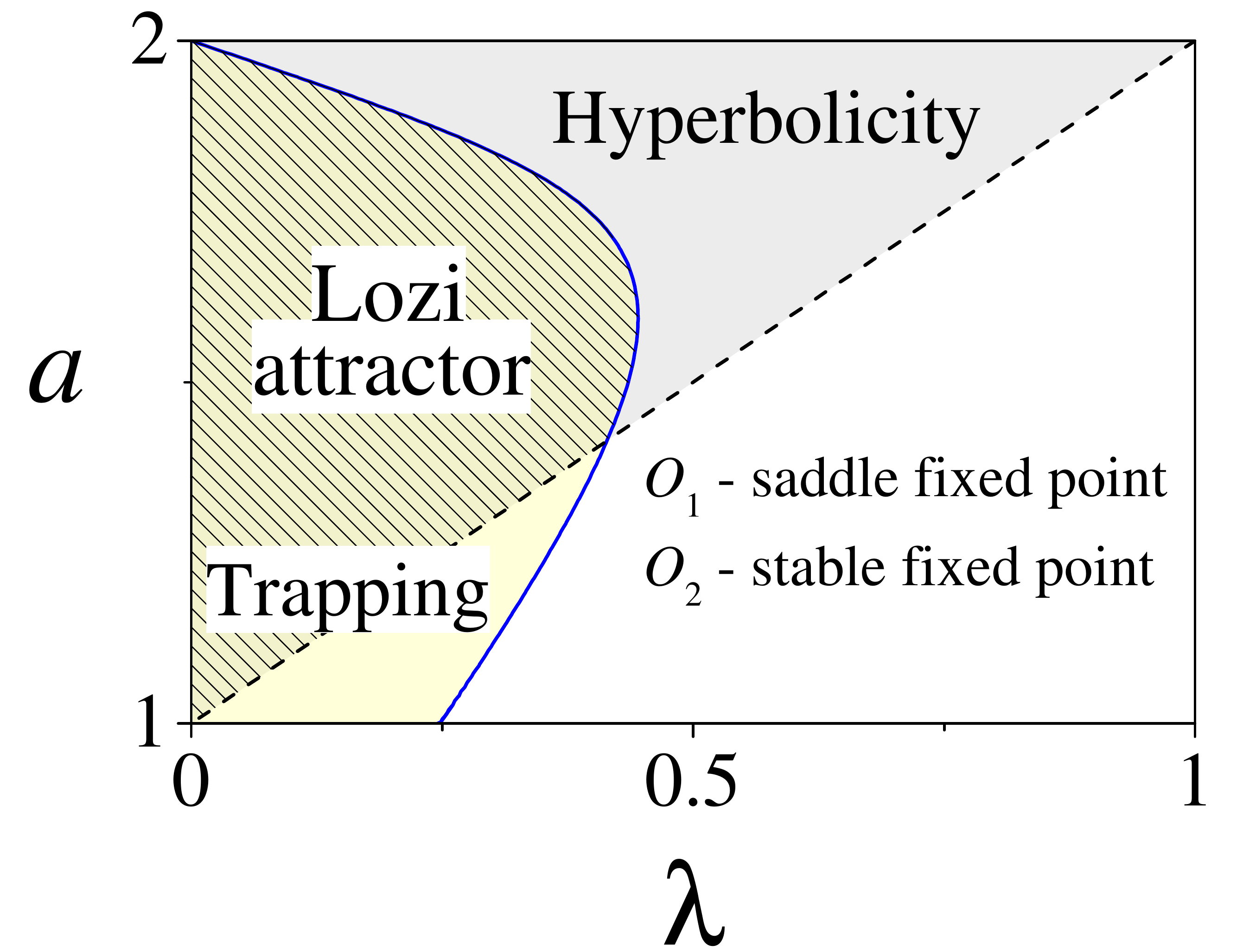}
    \caption{Parameters' plane $(\lambda, d)$ of the Lozi map~\eqref{sys:map_f}, \eqref{eq:g(x)_Lozi}.
    Parameters from yellow region corresponds to existence of invariant domain $G_t: fG_t\subset G_t$. This region is bounded by bifurcation of homoclinic orbit to saddle $O_1$ (blue curve), defined by equation $H(\lambda, d)=0$ in \eqref{eq:Lozi_homoclinic}.  Corresponding phase picture see in Fig.~\ref{fig:Lozi_homoclinic}.
    At black dashed line $a=1+\lambda$ fixed point $O_2$ becomes saddle (see conditions~\eqref{eq:conditions_d>0}).
    Grey region corresponds to hyperbolic phase space of the map.
    Intersection of yellow and gray regions (dashed area) corresponds to the existence of strange Lozi attractor.}
    \label{fig:Lozi_parameter_plane}
\end{figure}

Note that results of this sections are well-known (see \cite{Misiurewicz1980} and refs. therein). Here we present simple proof of the next statement.
\begin{theorem}
Let parameters $\lambda$ and $a$ satisfy the next two conditions
\begin{equation}
    a>1+\lambda,\quad \lambda<1,
    \label{eq:a>1_Lozi}
\end{equation}
and \begin{equation}
H(\lambda,a)\triangleq a(2\lambda-a+2)\sqrt{a^2-4\lambda}
+a(2\lambda^2-6\lambda-a^2+2a)+4\lambda^2\geq 0.
\label{eq:Lozi_homoclinic}
\end{equation}
Then the Lozi map~\eqref{sys:map_f}, \eqref{eq:g(x)_Lozi} has singular hyperbolic Lozi attractor (see Fig.~\ref{fig:Lozi_parameter_plane}). 
\end{theorem}

\textit{Proof.}
The derivative of function \eqref{eq:g(x)_Lozi}
\begin{equation}
    d(x)=-a\,\textrm{sign}\,x-1-\lambda
    \label{eq:d(x)_Lozi}
\end{equation}
satisfies conditions of hyperbolicity \eqref{eq:theorem_cond_d} in one and the same parameter domain \eqref{eq:a>1_Lozi} for $G^-(x>0)$ and $G^+(x<0)$.

Note that in this case in each point of $G$ expanding cone $K^u$ \eqref{eq:cones} contains two constant eigenvectors $\textbf{V}_1\big\vert_{x>0}$ and $\textbf{V}_2\big\vert_{x<0}$, and contracting cone $K^s$ contains two constant eigenvectors $\textbf{V}_1\big\vert_{x<0}$ and
$\textbf{V}_2\big\vert_{x>0}$.

Thereafter to prove the existence of Lozi attractor it is sufficient to find the part of parameter domain \eqref{eq:a>1_Lozi} corresponding to existence of trapping region $G_t$, $fG_t\subset G_t$.

Map~\eqref{sys:map_f},~\eqref{eq:g(x)_Lozi} has two saddle type fixed points $O_1\left(\frac{1}{1+\lambda-a},0\right)$ and $O_2\left(\frac{1}{1+\lambda+a},0\right)$ with multipliers $\mu_{1,2}$ 
\begin{equation*}
    \mu_{1,2}=\frac{1}{2}\left(
    -a\,\textrm{sign}\,x\pm\sqrt{a^2-4\lambda}
    \right),
\end{equation*}
satisfying conditions \eqref{eq:conditions_d>0}.
Eigenvectors $\textbf{V}_{1,2}=(1,\alpha_{1,2}(x))$, where coordinates \eqref{eq:multypliers_vectors} take the form
\begin{equation*}
    \alpha_{1,2}(x)=\frac{1}{2}(a\,\textrm{sign}\,x+2\lambda\pm\sqrt{a^2-4\lambda}).
\end{equation*}

As far as the eingenvectors are constant in each region, map~\eqref{sys:map_f},~\eqref{eq:g(x)_Lozi}
has unstable and stable invariant foliations in regions $x<0$ and $x>0$ parallel to separatrices $W_{1,2}$ of fixed points $O_1$ and $O_2$
\begin{equation}
\begin{array}{cl} 
W_{1,2}\big\vert_{x<0}: y=\alpha_{1,2}
\left(x+\frac{1}{a-1-\lambda} \right),
\\
W_{1,2}\big\vert_{x>0}: y=\alpha_{1,2}
\left(x-\frac{1}{a+1+\lambda} \right).
\end{array}
\label{eq:separ_Lozi}
\end{equation}

We construct continuation of unstable separatrix of point $O_1$ using \eqref{eq:separ_Lozi} as its first part. 
This part $W_1\vert_{x<0}$ intersects discontinuity line $\lbrace x=0, y\in\mathbb{R}^1 \rbrace$ in the point $M_1(0,\frac{\alpha^-_1}{a-1-\lambda})$, where $\alpha^-_1=\alpha_1\big\vert_{x<0}$.
Then continuation of $W_1\big\vert_{x<0}$ to region $x>0$ is line segment $M_1M_2$, where point $M_2=f\big\vert_{x<0}M_1=(x_2,y_2)$, $x_2=\frac{\alpha^-_1}{a-1-\lambda}+1$, $y_2=\lambda x_2$.

The next continuation of $W_1$ is line segment $M_2M_3$, where $M_3=f\big\vert_{x>0}(M_2)$, $M_3=(x_3,y_3)$, $x_3=y_2+1-(a+\lambda)x_2$, $y_3=\lambda(y_2+1-(a+1+\lambda)x_2)$.
For parameter region \eqref{eq:a>1_Lozi}, $M_3$ lies in region $x<0$ due to $x_3=\lambda x_2+1-(a+\lambda)x_2=1-a-\frac{a\alpha_1^-}{a-1-\lambda}<0$.

Introduce segment $I=\lbrace x=x_3, \alpha_2(x_3+\frac{1}{a-1-\lambda})\leq y \leq y_3 \rbrace$, connecting $M_3$ and unstable separatrix $W_2\big\vert_{x<0}$.

The area $G_t$ bounded by $W_1\big\vert_{x<0}\cup M_1M_2\cup M_2M_3\cup I$ is trapping region, $fG_t\subset G_t$, under condition \eqref{eq:Lozi_homoclinic}, which is redone inequality $\alpha_2(x_3+\frac{1}{a-1-\lambda})\leq y_3$.   

Hence, common parameter region \eqref{eq:a>1_Lozi},  \eqref{eq:Lozi_homoclinic} corresponds to Lozi attractor.
$\Box$

\begin{remark}
For $H(\lambda,a)=0$ point $M_3$ reaches stable separatrix $W_2\big\vert_{x<0}$ and becomes a point of homoclinic orbit. 
After this homoclinic bifurcation for $H(\lambda,a)<0$ point $M_3$ lies under $W_2\big\vert_{x<0}$, thus area $G_t$ is no longer trapping region.
\end{remark}

\begin{figure}
    \centering
    \includegraphics[width=0.75\textwidth]{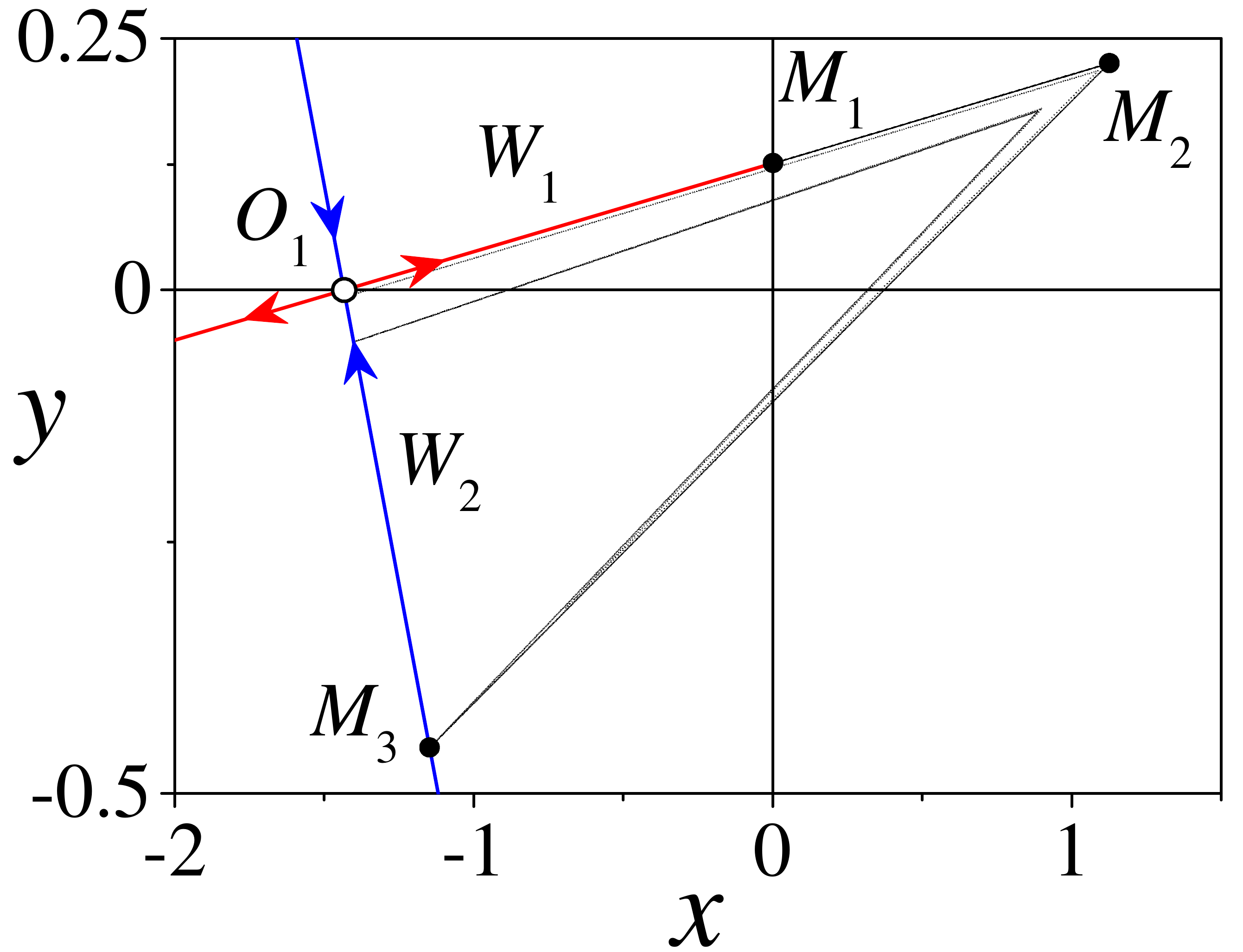}
    \caption{Homoclinic orbit in the Lozi map~\eqref{sys:map_f}, \eqref{eq:g(x)_Lozi} for $\lambda=0.2$ and $a=1.905$.
    Unstable and stable manifolds $W_{1,2}$ of saddle $O_1$ are red and blue lines, respectively. Black broken line is the continuation of $W_1$, which forms homoclinic orbit.}
    \label{fig:Lozi_homoclinic}
\end{figure}

Fig.~\ref{fig:Lozi_homoclinic} illustrates that breakdown of Lozi attractor occurs as a result of either hyperbolicity condition loss or via homoclinic bifurcation.

\section{\label{sec:Lozi_Hybrid}Attractor of nonlinear hybrid Lozi-H\'{e}non map}

Consider map (6) in the case of function
\begin{equation}
    U(x)=1-a(l|x|+(1-l)x^2),
    \label{eq:U(x)_Hybrid}
\end{equation}
which is the Lozi map for $l=1$ and H\'{e}non map for $l=0$.
As before we transform map (6) to map~\eqref{sys:map_f} with function
\begin{equation}
    g(x)=1-al|x|-(1+\lambda)x-a(1-l)x^2.
    \label{eq:g(x)_Hybrid}
\end{equation}
The condition of hyperbolicity \eqref{eq:theorem_cond_d} become the same for $G=G^-\cup G^+$ condition for $l>0$
\begin{equation}
    al>1+\lambda,
    \label{eq:al>1_Hybrid}
\end{equation}
and condition for H\'{e}non map $l=0$
\begin{equation}
    |x|>\frac{1+\lambda}{2a}.
    \label{eq:x>1_Hybrid}
\end{equation}

Note that due to \eqref{eq:al>1_Hybrid} when $l$ tends to zero when H\'{e}non part dominates in $g(x)$ parameter $a$ in \eqref{eq:al>1_Hybrid} increases to infinity.

Condition \eqref{eq:al>1_Hybrid} eliminates extremums of $g(x)$ with $g'(x)=0$ and force $g(x)$ to have maximum at point $(x=0, y=1)$ with discontinuous derivative.

Map~\eqref{sys:map_f}, \eqref{eq:g(x)_Hybrid}
has saddle points $O_{1,2}$ $(x_{1,2}, 0)$,
where 
\begin{equation}
    x_{1,2}^s=\frac{1}{2a(1-l)}\left[
    -1-\lambda\pm\left(al-\sqrt{(al\mp(1+\lambda))^2+4a(1-l)}
    \right)
    \right].
    \label{eq:x_12^s_Hybrid}
\end{equation}

Here separatrices of saddle $O_1$ are non-linear in initial region $x<0$ and defined by nonlinear functions
\begin{equation}
    W_{1,2}\big\vert_{G^-}:\;y=\varphi_{1,2}(x-x_1^s).
\end{equation}

\begin{theorem}
There exists a value $l^*\in(0,1)$ such that for $l^*<l<1$ hybrid Lozi-H\'{e}non map has strange hyperbolic attractor if condition \eqref{eq:al>1_Hybrid} and inequality
    \begin{equation}
        -y_3\leq (1-\lambda)(x_3-x_1^s)
        \label{eq:th_Hybrid}
    \end{equation}
    are valid. 
    \label{th:attr_Hybrid}
\end{theorem}

\textit{Proof.}
Let condition of hyperebolicity \eqref{eq:al>1_Hybrid} is valid.
Consider sufficient conditions for existence of trapping region.

In parameter domain \eqref{eq:al>1_Hybrid} cones \eqref{eq:cones} give the next bounds of separatrices
\begin{equation}
    \begin{array}{c}
         0<\varphi_{1}(x-x_1^s)<\lambda(x-x_1^s),\\
         \varphi_{2}(x-x_1^s)<-(1-\lambda)(x-x_1^s).
    \end{array}
    \label{eq:0<phi_1_Hybrid}
\end{equation}

Hence, unstable separatrix $W_1\big\vert_{x<0}$ intersects the line $\lbrace x=0, y\in\mathbb{R}^1 \rbrace$ in point $M_1(0, y_1)$, where
\begin{equation}
    0<y_1<\bar{y}_1=-\lambda x_1^s.
\end{equation}

The image $f\big\vert_{x<0}(M_1)=M_2(x_2,y_2)$ has the coordinates
\begin{equation}
    x_2=y_1+1,\quad y_2=\lambda(y_1+1),
\end{equation}
and the first piece of continuation of separatrix $W_1\big\vert_{x<0}$ to region $x>0$ is curve $\wideparen{M_1M_2}$ connecting $M_1$ and $M_2$.
The next piece is curve $M_2M_3$ connecting $M_2$ with its image $M_3=f\big\vert_{x>0}(M_2)$, $M_3(x_3,y_3)$, where 
\begin{equation}
    \begin{array}{l}
         x_3=1-al(y_1+1)-a(1-l)(y_1+1)^2, \\
         y_3=\lambda-\lambda(al+1)(y_1+1)-\lambda a(1-l)(y_1+1)^2, 
    \end{array}
    \label{eq:x_3_Hybrid}
\end{equation}

From \eqref{eq:x_3_Hybrid} we obtain the next inequalities
\begin{equation}
\begin{array}{c}
    x_3<x_3\big\vert_{y_1=0}=1-a<0,\\
    y_3<y_3\big\vert_{y_1=0}=-\lambda a<0.
\end{array}
\label{eq:x_3<x_3_Hybrid}
\end{equation}
implying that $M_3\in G^-(x<0)$ and the second piece of separatrix continuation $\wideparen{M_2M_3}$ intersects region $x<0$.

Similarly to the case of Lozi map $l=1$ we introduce segment $I_l=\left\lbrace x=x_s, \varphi_2(x_3-x_1^s)<y<y_3 \right\rbrace$.
As image $f\big\vert_{x<0}I_l\in\left\lbrace x_3<0 \right\rbrace$, so area $G_l$ bounded by closed curve $\partial G_l=W_1\big\vert_{x<0}\cup\wideparen{M_1M_2}\cup\wideparen{M_2M_3}\cup I_l$ is trapping region of map~\eqref{sys:map_f},~\eqref{eq:g(x)_Hybrid}.

Using formulas in \eqref{eq:x_12^s_Hybrid}-\eqref{eq:x_3_Hybrid}
we obtain that the region of parameters \eqref{eq:th_Hybrid} is not empty and exists only for $0<a<a^*(l,\lambda)$.
Then the condition of hyperbolicity \eqref{eq:al>1_Hybrid} and condition of trapping domain existence \eqref{eq:th_Hybrid} are compatible for some value $l^*$ defined by equation $a^*(\lambda,l)l=1+\lambda$.
$\Box$

\begin{remark}
Obviously Theorem~\ref{th:Hyper} does not work in the case of H\'{e}non map for $l=0$ due to breach of hyperbolicity \eqref{eq:al>1_Hybrid} for the whole domain $G=G^-\cup G^+$.
However condition \eqref{eq:th_Hybrid} for trapping domain is valid for $l=0$, and condition of hyperbolicity \eqref{eq:x>1_Hybrid} can be used for nonwandering subset fully located in domain \eqref{eq:x>1_Hybrid}.
\end{remark}

\section{\label{sec:Belykh} Discontinuous function and Belykh attractor}
Belykh map \cite{BelykhMapScolarpedia} and its multidimensional generalizations \cite{PesinDynamicalSys1992,Sataev1992,Sataev1999}. were considered in terms of ergodic theory in order to prove existence of continuous invariant measure \cite{BelykhMapScolarpedia,AfraimovichChernovSataev1995,PesinDynamicalSys1992,Sataev1992,Sataev1999,SchmelingErgoticTheory1998,SchmelingTroubetzkoy1998}. Here we consider dynamical properties of Belykh attractor and present its new features for the map $f$ of cylinder $S^1\times R^1$.

We introduce sawtooth 2-periodic function comming from PLL modeling
\begin{equation}
g(x)=a(x-(2k-1))\quad\textrm{for}\quad x \in (2k-2, 2k),\quad k \in \mathbb{Z}.
\label{pll:g_first}
\end{equation}

In the case $x\in\mathbb{S}^1$ we can rewrite $g(x)$ in the form
\begin{equation*}
g(x)=a(x-1),\quad x(\textrm{mod}\;2).
\end{equation*}

First we consider the map~\eqref{sys:map_f} at two periods $(k=0, 1)$ of function $g(x)$ which takes the form
\begin{equation}
g(x)=
\left\{
\begin{array}{cl}
a(x+1),\quad x<0,\\
\\
a(x-1),\quad x\geq 0.
\end{array}
\right.
\label{pll:g_new}
\end{equation}
Here for convenience we consider $x \in\mathbb{R}^1$ instead $x \in [-2, 2]$.

Matrix $A$ in \eqref{sys:map_f_L} having the form $A=\begin{pmatrix} 1+a & 1\\ \lambda a & \lambda \end{pmatrix}$ is constant for any $x \in\mathbb{R}^1$, and as $d(x)\equiv a >0$ and trace $s\triangleq 1+a+\lambda$ of matrix $A$, so due to \eqref{eq:multypliers_vectors} eigenvalues $\mu_1 >1$, $0<\mu_2<1$ and eigenvectors $\textbf{V}_{1,2}=(1, \alpha_{1,2})$ coordinates $\alpha_{1,2}$ satisfy conditions
\begin{equation}
0<\alpha_1<1, \quad \alpha_2<-(1+\lambda)<0.
\end{equation}

As far as constant eigenvectors $\textbf{V}_{1,2}$ are the same for any $x$ map $f$ has unstable $\Phi_c^u$ and stable $\Phi_c^s$ invariant foliations
\begin{equation}
\begin{array}{cl}
\Phi_c^u=\{x, y\,|\, y=\alpha_1(x+c_1)\},\\
\\
\Phi_c^s=\{x, y\,|\, y=\alpha_2(x+c_2)\},
\end{array}
\label{pll:inv_foliations}
\end{equation}
where constants $c_{1,2}\in\mathbb{R}^1$.
Note that cones $K^u$ and $K^s$ 
are superfluous in the case of map~\eqref{sys:map_f},~\eqref{pll:g_new}
because each of them contain one constant vector everywhere.

\begin{figure}
    \centering
    \includegraphics[width=0.75\textwidth]{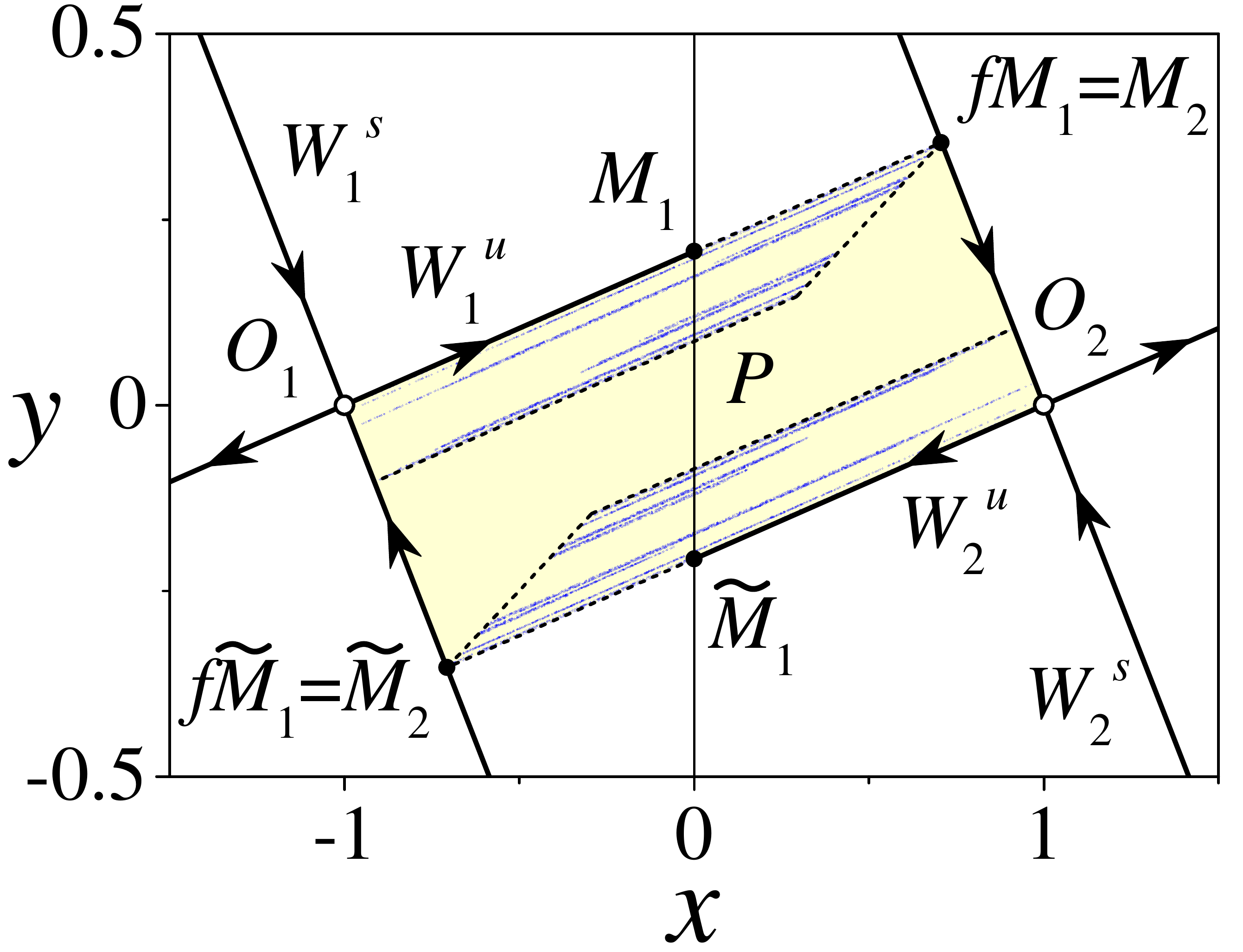}
    \caption{Phase portrait of Belykh map~\eqref{sys:map_f}, \eqref{pll:g_new} for bifurcational case $\lambda=a=0.5$ (see condition \eqref{eq:a<1}) when point $M_1$ maps to stable manifold $W_2^s$, $fM_1=M_2$. This is thr edge case when parallelogram $P$ (yellow) forms trapping region $fP\subset P$.
    This parallelogram maps into two trapezoids (areas bounded with dashed lines): upper $fP\big\vert_{x<0}$ and bottom $fP\big\vert_{x>0}$. Blue points are singular hyperbolic Belykh attractor.}
    \label{fig:Belykh_map}
\end{figure}

The map $f$ \eqref{sys:map_f}, \eqref{pll:g_new} has two similar fixed points $O_{1,2}(\pm 1, 0)$ of saddle type with unstable $W_{1,2}^u$ and stable $W_{1,2}^s$ separatrices (see Figure~\ref{fig:Belykh_map})
\begin{equation}
\begin{array}{cl}
W_{1,2}^u=\{x, y\,|\, y=\alpha_1(x\pm 1)\},\\
\\
W_{1,2}^s=\{x, y\,|\, y=\alpha_2(x\pm 1)\},
\end{array}
\label{pll:separatrices}
\end{equation}
lying in invariant foliations~\eqref{pll:inv_foliations}.

Denote points $M_1(x_1,y_1)=W_1^u\cap \{x=0\}$, where $x_1=0, y_1=\alpha_1$ and $M_2(x_2,y_2)=W_2^s\cap\left\lbrace x>0, y=\alpha_1(x+1)\right\rbrace$, where 
\begin{equation*}
x_2=\frac{s-2\lambda}{\sqrt{s^2-4\lambda}},\quad y_2=\alpha_1(x_2+1).
\end{equation*}

\begin{theorem}
    For
    \begin{equation}
     a\leq 1-\lambda
     \label{eq:a<1}
    \end{equation}
map~\eqref{sys:map_f}, \eqref{pll:g_new} has singular hyperbolic Belykh attractor.
\label{pll:th:Belykh_attr}
\end{theorem}

\textit{Proof.} 
The image $f\big\vert_{x<0}(M_1)=\big(\bar{x}_1, \bar{y}_1)=(\alpha_1+a, \lambda(\alpha_1+a)\big)$, does not reach point $M_2$ when $\bar{x}_1<x_2$, i.e. for
\begin{equation}
    \alpha_2+a<\frac{s-2\lambda}{\sqrt{s^2-4\lambda}}.
    \label{eq:alpha_2<a}
\end{equation}
Symmetrically $f\big\vert_{x>0}(\bar{M}_1)$ does not reach point $\bar{M}_2$ for \eqref{eq:alpha_2<a}.
Inequality~\eqref{eq:alpha_2<a} takes the form
\begin{equation}
    (2-s)(s+\sqrt{s^2-4\lambda})>0,
    \label{eq:2-A}
\end{equation}
which is valid under condition \eqref{eq:a<1}.
Hence parallelogram $P=O_1M_2O_2\bar{M}_2$ is invariant, $fP\subset P$, and there exists an attractor $\mathcal{A}\subset P$.

Folliations \eqref{pll:inv_foliations} cross domain $P$ for $c_{1,2}\in [-1,1]$ with stretch $\mu_1>1$ along $\Phi_c^u$ and shrink $\mu_2<1$ along $\Phi_c^s$.
Hence each orbit of attractor $\mathcal{A}$ is uniformly hyperbolic with Lyapunov exponents $h_1=\ln \mu_1>0$ and $h_2=\ln \mu_2<0$.
Singularity of attractor $\mathcal{A}$ is due to the fact that a set of parameters $(a,\lambda)$ corresponding to $\mathcal{A}\cap\lbrace x=0 \rbrace\neq\emptyset$ is a set of ``border collision'' bifurcations.
For example, period $p$ orbit disappears when its point reaches discontinuity line $\lbrace x=0 \rbrace$ under a parameter change.
$\Box$

\begin{remark}
In basis of eigenvectors map~\eqref{sys:map_f},~\eqref{pll:g_new} can be written in the form
\begin{equation}
    \begin{array}{lll}
         \bar{u}_1=\mu_1(u_1+1)+1,& 
         \bar{u}_2=\mu_2(u_2+1)+1,
         &v<0,\\
         \\
         \bar{u}_1=\mu_1(u_1-1)-1,& 
         \bar{u}_2=\mu_2(u_2-1)-1,
         &v>0,
    \end{array}
    \label{sys:Belykh_map_in_basis}
\end{equation}
where $v=u_1+k u_2$, $0<k<1$.
In the literature \cite{AfraimovichChernovSataev1995}, Belykh map is met just in the form~\eqref{sys:Belykh_map_in_basis}.
\end{remark}

\subsection{Periodic case}

Consider the case of 2-periodic function $g(x)$ \eqref{pll:g_first} which becomes bounded function.

\begin{lemma}
Annulus $R=\lbrace x,y\;\vert\; x\in \mathbb{S}^1, |y|<\frac{a}{1-\lambda}\rbrace$
is absorbing domain of map $f$~\eqref{sys:map_f}, \eqref{pll:g_first}.
\label{th:lemma_annulus_periodic}
\end{lemma}

\textit{Proof.}
Function $g(x)=a(x-1)$ is bounded at $\mathbb{S}^1$: $|g(x)|<a$ at the period $0<x<2$.
Then $|\bar{y}|<\lambda(|y|+a)$ and variable $|y|$ decreases at each iterate for $\lambda(|y|+a)<|y|$, i.e. for $|y|\geq \frac{a}{1-\lambda}$. The latter implies that annulus $R$ is absorbing domain.$\Box$

\begin{theorem}
    For $a>1-\lambda$, $0<\lambda<1$ 
    map~\eqref{sys:map_f}, \eqref{pll:g_first}
    has singular hyperbolic attractor consisting of three components: oscillating $\mathcal{A}_o\subset P$, rotating $\mathcal{A}_r\subset R\setminus P$ and rotating with oscillations $\mathcal{A}_{ro}\subset R$.
    \label{th:Belykh_map_periodic}
\end{theorem}

\textit{Proof.}
From Lemma~\ref{th:lemma_annulus_periodic} it follows that $fR\subset R$ and the map has an attractor $\mathcal{A}$.
Due to existence of folliations $\Phi_c^u$ and $\Phi_c^s$ the attractor is singular hyperbolic. 
For $a>1-\lambda$ image $fM_1$ overjumps point $M_2$ and $fP\setminus P\neq\emptyset$ is triangle $\Delta_1=M_2fM_1M_4$ (or trapezoid) for $y>0$
and symmetrical $\bar{\Delta}_1$ for $y<0$.
Image $f^k\Delta$ beginning from some $k$ starts to enwind the cylinder forming rotating attractor $\mathcal{A}_r$ and intersecting the triangle $\Delta_2=M_1M_3M_2$.
Obviously $f^k\Delta_2\subset P$, where $k\geq 1$,
i.e. trajectories from $\Delta_2$ reach region $P$.
Almost all trajectories $f^k\Delta_2$, $k>k^*$,  leave region $P$ and form attractor $\mathcal{A}_{ro}$. $\Box$

\begin{figure}
    \centering
    \includegraphics[width=0.75\textwidth]{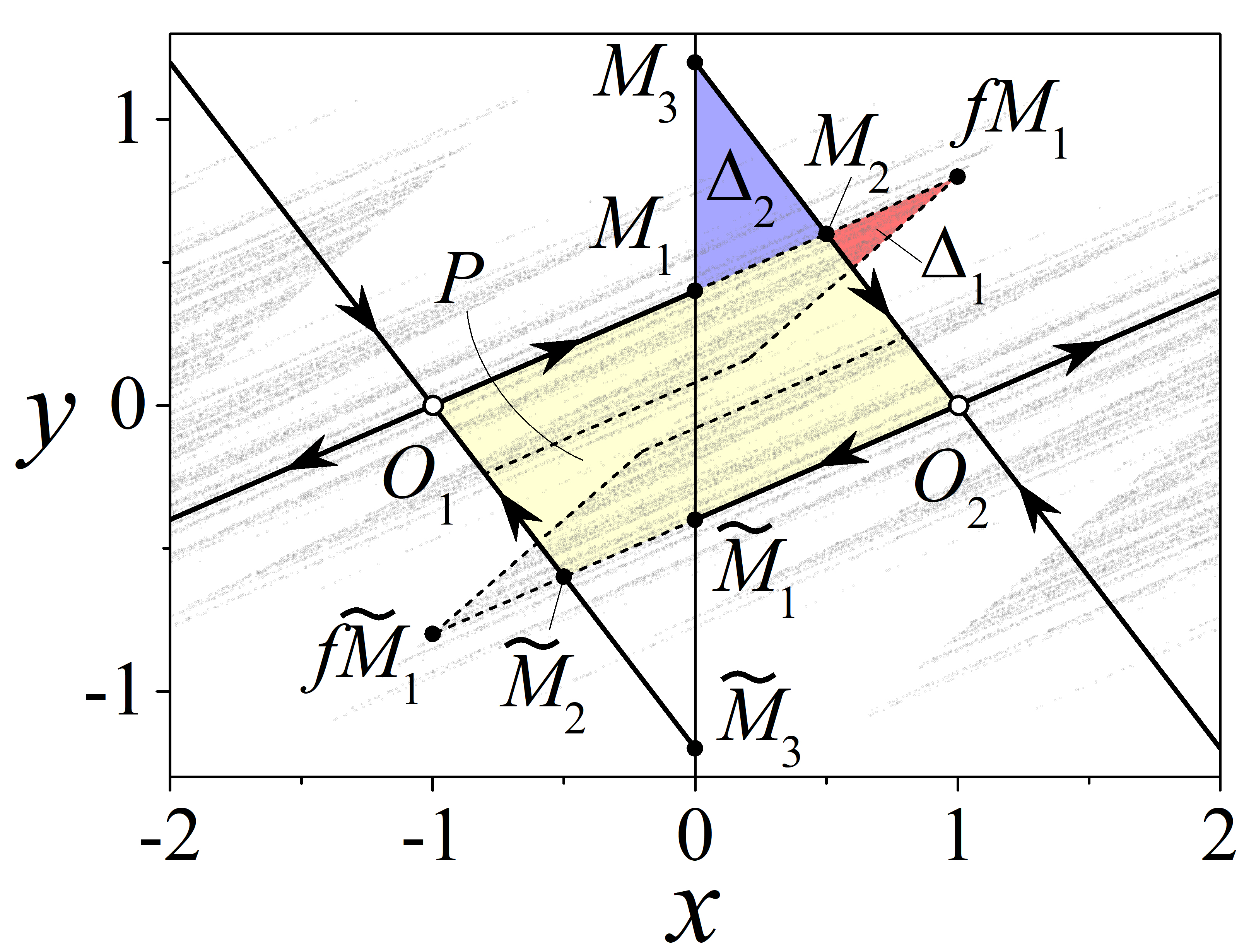}
    \caption{Phase portrait of Belykh map~\eqref{sys:map_f} for periodic function $g(x)$ \eqref{pll:g_first} under conditions of Theorem~\ref{th:Belykh_map_periodic} (with parameters $\lambda=0.8$, $a=0.6$).
    Attractors $\mathcal{A}_{o}$, $\mathcal{A}_r$ and $\mathcal{A}_{ro}$ are light gray dots.
    Parallelogram $P$ (yellow) maps into two trapezoids (dashed contours): upper $fP\big\vert_{x<0}$ and bottom $fP\big\vert_{x>0}$.
    Red triangle $\Delta_1$ plays the role of ``output gate'' from $P$ for trajectories of $A_ro$, and blue triangle $\Delta_2$ plays the role of ``input gate'' to $P$.
    Coordinates of points: $M_1(0, 0.4)$, $M_2(0.5, 0.6)$, $M_3(0, 1.2)$ $fM_1(1, 0.8)$, $\tilde{M}_1(0, -0.4)$, $\tilde{M}_2(-0.5, -0.6)$, $f\tilde{M}_1(-1, -0.8)$. Points $\tilde{M}$ are odd symmetric to points $M$.}
    \label{fig:attr_rotating}
\end{figure}

\textbf{Corollary.} 
For $a>0$, $0<\lambda<1$ map $f$ has singular hyperbolic attractor.
Indeed for $a\leq 1-\lambda$ that is Belykh attractor and for $a>1-\lambda$ three attractors $\mathcal{A}_o$, $\mathcal{A}_r$ and $\mathcal{A}_{ro}$ from Theorem~\ref{pll:th:Belykh_attr}.

\begin{remark}
Attractor $\mathcal{A}_o$ is a Cantor set similar to that of Smale's horseshoe.
Unlike attractors $\mathcal{A}_r$ and $\mathcal{A}_{ro}$ it's basin haz zero Lebesgue measure. 
\end{remark}

\begin{remark}
The model of continuous time PLL~\eqref{sys:PLL} has simple dynamics.
Indeed, Lyapunov function
\begin{equation*}
    V=\frac{u^2}{2}-\omega_2\int\limits_0^xg(\xi)d\xi,
\end{equation*}
for system~\eqref{sys:PLL},~\eqref{pll:g_first} with derivative
\begin{equation*}
    \dot{V}=-\alpha u^2-\omega_1\omega_2g^2(x)<0\quad\textrm{for}\; x\neq 2k-1, u\neq 0,
\end{equation*}

Hence system~\eqref{sys:PLL} is globally asymptotically stable so that all trajectories, besides stable separatrices, are attracted by the point $(x=u=0)$ - glued stable equilibrium (see Figure~\ref{fig:PLL_phase}).

Under transition from system~\eqref{sys:PLL} to map~\eqref{sys:map_f},~\eqref{pll:g_first} this zero point immediately becomes strange attractor which ``amplitude'' is of order $\alpha_1 +a \backsim h$ (see \eqref{eq:change_for_Belykn_map}).
\begin{figure}[h!]
    \centering
    \includegraphics[width=0.75\textwidth]{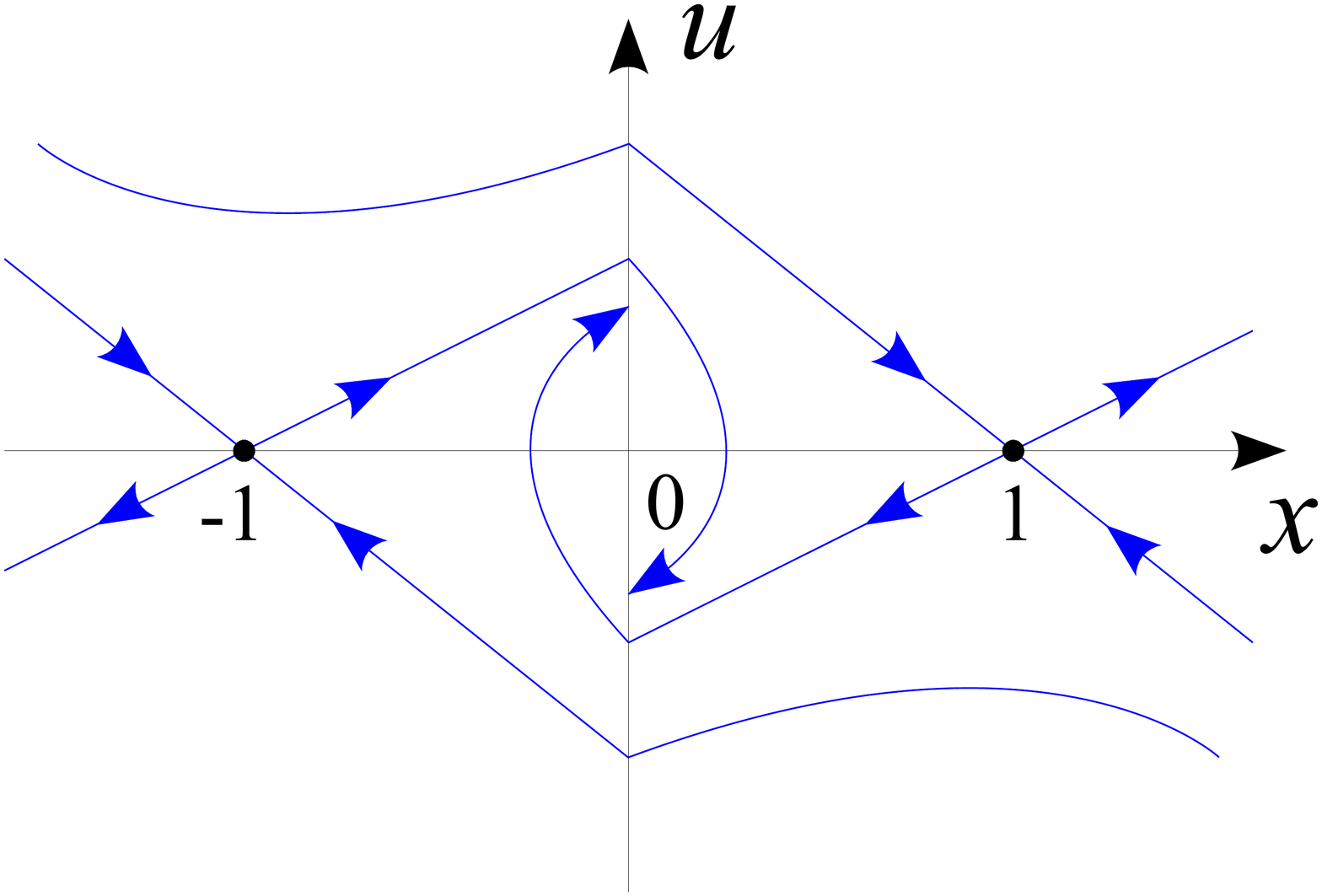}
    \caption{Phase portrait of phase-locked loop system~\eqref{sys:PLL} for $\nu$=0.}
    \label{fig:PLL_phase}
\end{figure}
\end{remark}

\section{Acknowledgements}
This work was supported by the Russian Science Foundation under grant No. 22-21-00553.

\end{document}